\include{psfig.sty}
\documentclass[11pt,a4paper]{article}
\usepackage{graphicx}
\usepackage{amsfonts}
\usepackage{latexsym}
\newcommand{\beql}[1]{\begin{equation}\label{#1}}
\newcommand{\eeq}{\end{equation}}
\newcommand{\comment}[1]{}

\newcommand{\eqref}[1]{{\rm (\ref{#1})}}

\newcommand{\Abs}[1]{{\left|{#1}\right|}}

\newcommand{\Qed}{\ \\\mbox{$\Box$}}

\newcounter{open}
\newcounter{dfn}
\def\thedfn{\arabic{dfn}}

\newcounter{obs}
\def\theobs{\arabic{obs}}

\newcounter{thm}
\newcounter{mysec}
\newcounter{mysubsec}[mysec]
\newtheorem{theorem}{Theorem}

\begin{document}

\begin{center}
{\Large \bf A remark on perturbations of sine and cosine sums}\\
\ \\
{\sc Mihail N. Kolountzakis}\\
\ \\
\small December 1999
\end{center}

\noindent
Consider a collection $\lambda_1<\cdots<\lambda_N$ of distinct
positive integers and the quantities
$$
M_1 = M_1(\lambda_1,\ldots,\lambda_N) = 
 \max_{0\le x \le 2\pi} \Abs{\sum_{j=1}^N \sin{\lambda_j x}}
$$
and
$$
M_2 = M_2(\lambda_1,\ldots,\lambda_N) = - \min_{0\le x \le 2\pi}
 \sum_{j=1}^N \cos{\lambda_j x}.
$$
One is interested in frequencies $\lambda_j$ which make the quantities
$M_1$ and $M_2$ small as $N \to \infty$.
It is trivial that $M_1 \ge c N^{1/2}$ but it is much harder
to show even that $M_2 \to \infty$.
(In this note $c$ denotes an absolute positive constant,
not necessarily the same in all its occurences.)

It is a result of Bourgain \cite{Bourgain,Kahane}
that $M_1$ may become $O(N^{2/3})$ and it is very easy to construct
a collection $\lambda_j$ which gives $M_2 = O(N^{1/2})$, which is the
conjectured optimal. For minimizing $M_2$ it is also possible to
have the collection of frequencies relatively well packed, that
is with $\lambda_N \le 2N$ \cite{Kolountzakis}, while for
any $\epsilon>0$ and for
any collection $\lambda_j$ that makes $M_1 = O(N^{1-\epsilon})$ 
one can easily see that $\lambda_N$ is super-polynomial in $N$.

Prompted by a discussion with G. Benke we prove that collections of
frequencies $\lambda_j$ which have $M_1 = o(N)$ or $M_2 = o(N)$ are unstable,
in the sense that one can perturb the $\lambda_j$ by one each and
get $M_1 \ge c N$ and $M_2 \ge c N$.
\begin{theorem}\label{th:main}
{\rm (i)} Suppose that $M_1(\lambda_1,\ldots,\lambda_N) = o(N)$.
Then there exists a choice of $\epsilon_j = \pm 1$, $j=1,\ldots,N$,
so that $M_1(\lambda_1+\epsilon_1,\ldots,\lambda_N+\epsilon_N) \ge cN$.\\
{\rm (ii)} Suppose, similarly, that $M_2(\lambda_1,\ldots,\lambda_N) = o(N)$.
Then there exists a choice of $\epsilon_j = \pm 1$, $j=1,\ldots,N$,
so that $M_2(\lambda_1+\epsilon_1,\ldots,\lambda_N+\epsilon_N) \ge cN$.
\end{theorem}
{\bf Remark.}
It is not always the case that the perturbed frequencies are all
distinct but is frequently so and, in any case, at most two
may overlap at any given integer.

\noindent
{\bf Proof.}
(i)\ \  Write
\begin{eqnarray*}
\sum_{j=1}^N \sin{(\lambda_j+\epsilon_j)x} &=&
 \sum_{j=1}^N \sin{\lambda_j x} \cos{\epsilon_j x} +
 \sum_{j=1}^N \cos{\lambda_j x} \sin{\epsilon_j x}\\
 &=& \cos x \sum_{j=1}^N \sin{\lambda_j x} +
     \sin x \sum_{j=1}^N \epsilon_j \cos{\lambda_j x}\\
 &=& {\rm I} + {\rm II}.
\end{eqnarray*}
We have ${\rm I} = o(N)$.
For $\lambda>10$, say, we have
$$
{4\over \pi} \int_{\pi/4}^{\pi/2} \Abs{\cos{\lambda x}}~dx \ge c.
$$
From this we deduce that
$$
{4\over \pi} \int_{\pi/4}^{\pi/2} \sum_{j=1}^N \Abs{\cos{\lambda_j x}}~dx
	\ge cN,
$$
hence there exists $x_0 \in [{\pi\over4},{\pi\over2}]$ such that
$$
\sum_{j=1}^N \Abs{\cos{\lambda_j x_0}} \ge cN.
$$
Choose then $\epsilon_j = {\rm sgn}\,({\cos{\lambda_j x_0}})$ to get
$$
\sum_{j=1}^N \epsilon_j \cos{\lambda_j x_0} \ge cN.
$$
Since $\sin x \ge 2^{-1/2}$ in $[{\pi\over4},{\pi\over2}]$
we get that ${\rm II} \ge cN$ at $x_0$, which gives the
required $M_1(\lambda_1+\epsilon_1,\ldots,\lambda_N+\epsilon_N) \ge cN$,
since ${\rm I} = o(N)$ everywhere.

(ii)\ \ The proof is similar. We write
\begin{eqnarray*}
\sum_{j=1}^N \cos{(\lambda_j+\epsilon_j)x} & = &
 \cos x \sum_{j=1}^N \cos{\lambda_j x} -
 \sin x \sum_{j=1}^N \epsilon_j \sin{\lambda_j x}\\
& = & {\rm I} - {\rm II}.
\end{eqnarray*}
For $x \in [{4\pi\over6},{5\pi\over6}]$ we have ${\rm I} \le o(N)$ so it is
enough to show that for some $x_0$ in the same interval ${\rm II}\ge cN$. 
To do that observe, as before, that
\beql{abs-mean}
{6\over\pi}\int_{4\pi/6}^{5\pi/6}\sum_{j=1}^N \Abs{\sin{\lambda_j x}} \ge cN,
\eeq
and choose $\epsilon_j = {\rm sgn}\,(\sin{\lambda_j x_0})$, where
$x_0$ is a point that makes the left hand side of \eqref{abs-mean}
large.
This implies the existence of $x_0 \in [{4\pi\over6},{5\pi\over6}]$ such that
${\rm II} \ge cN$ at $x_0$.
\Qed

Take a Bourgain's collection of frequencies $\lambda_j$ for which
$M_1 = O(N^{2/3})$.
Bourgain gave a randomized construction for these which produces
$N$ frequencies in the interval $[1,e^{N^{1/3}}]$.
By taking the inner product of the above sine sum with the conjugate
Dirichlet kernel $D_M^* = \sum_{j=1}^M \sin {jx}$,
whose $L^1$ norm is $c \log M$, we obtain
that the number of frequencies in the interval $[1,M]$ is at most
$c N^{2/3} \log M$, hence Bourgain's method does not go any further
than it has to (in the size of the $\lambda_j$).
One can even modify his randomized construction to produce
$\lambda_j \sim e^{j^{1/3}}$.
Hence the question becomes reasonable of whether it is just
the growth of the frequencies that achieves the result.
Our Theorem \ref{th:main} answers this in the negative.

\noindent
\ \\

\noindent
{\sc\small
Department of Mathematics, University of Crete, Knossos Ave.,
714 09 Iraklio, Greece.\\
E-mail: {\tt kolount@math.uch.gr}
}


{\bf Bibliography}

\begin{thebibliography}{99}
\vspace{-1.6cm}

\bibitem{Bourgain} J. Bourgain, Sur les sommes de sinus,
S\'em.\ Anal.\ Harm., Publ.\ Math.\ d'Orsay 84-01 (1984), exp.\ no 3.

\bibitem{Kahane} J.-P. Kahane, {\em Some random series of functions},
Cambridge Studies in Advanced Mathematics 5, 1985, Second Edition.

\bibitem{Kolountzakis} M.N. Kolountzakis, On nonnegative cosine polynomials
with nonnegative, integral coefficients,
Proc.\ Amer.\ Math.\ Soc.\ {\bf 120} (1994), vol.\ 1, 157-163.

\end{thebibliography}
\end{document}